\documentclass[reqno, draft]{amsart}
\usepackage{amssymb, amscd, amsfonts}
\usepackage[mathscr]{eucal}
\usepackage{graphicx}
\usepackage{amscd}

\newtheorem{theorem}{Theorem}[section]
\newtheorem{lemma}[theorem]{Lemma}
\newtheorem{corollary}[theorem]{Corollary}
\newtheorem{proposition}[theorem]{Proposition}

\theoremstyle{definition}
\newtheorem{definition}[theorem]{Definition}
\newtheorem{remark}[theorem]{Remark}

\DeclareMathOperator{\conv}{\rm conv}
\DeclareMathOperator{\den}{\rm den}
\DeclareMathOperator{\aff}{\rm aff}
\DeclareMathOperator{\cl}{\rm cl}
\DeclareMathOperator{\ext}{\rm ext}

\newcommand{\restrict}{\,{\mathbin{\vert\mkern-0.3mu\grave{}}}\,}

\newcommand{\remove}[1]{}

\title[Measure   theory in the geometry of
$GL(n,\mathbb Z) \ltimes \mathbb Z^{n}$]
{Measure   theory in the geometry of
$GL(n,\mathbb Z) \ltimes \mathbb Z^{n}$}

\author[D. Mundici]{Daniele Mundici$^\dag$}
\address[D. Mundici]{Department of Mathematics ``Ulisse Dini'' \\
University of Florence \\
viale Morgagni 67/A \\
50134 Florence \\
Italy}

\email{mundici@math.unifi.it }



 \date{\today}

\begin{document}

\begin{abstract}  The $n$-dimensional affine group over the integers
 is the group    $\mathcal G_n$
  of all affinities on $\mathbb R^{n}$ which leave the lattice $ \mathbb Z^{n}$
  invariant.  $\mathcal G_n$ yields a geometry in 
  the classical sense of the Erlangen Program.
In this paper we construct a $\mathcal G_n$-invariant measure
  on rational polyhedra in $\mathbb R^n$,   i.e., finite unions of simplexes
  with rational vertices  in $\mathbb R^n$,   
and prove its uniqueness.
Our main tool is given by 
the Morelli-W{\l}odarczyk factorization of birational toric maps in
blow-ups and blow-downs  (solution of the weak  
Oda  conjecture). 
\end{abstract}

\maketitle

\section{Introduction}
\label{section:introduction}
Following  Stallings \cite{sta}
and Tsujii \cite{tsu}, 
by a  {\it polyhedron}  $P$ in $\mathbb R^{n}  \,\, (n=1,2,\ldots)$
we mean the pointwise union of  a finite 
set of (always closed)  simplexes
$T_i$  in $\mathbb R^{n}$.
If the vertices of each  $T_i$ are in $\mathbb Q^n$,
$P$ is said to be  a
 {\it rational polyhedron. }
 $P$ need not be convex or connected: 
 as a matter of fact, every compact
 subset of $\mathbb R^n$ coincides with the intersection of all
 rational polyhedra in $\mathbb R^n$ containing it.
We will write
$
\mathscr P^{(n)}
$ for   the set of all 
rational polyhedra in $\mathbb R^n.$ 
 
We    let
$\mathcal G_n=GL(n,\mathbb Z) \ltimes \mathbb Z^{n}$ 
denote the group of transformations
of the form
$$
x\mapsto Ax + t \quad \mbox{for} \quad x\in \mathbb R^n,
$$
 where $t\in \mathbb Z^n$
 and  $A$ is an  $n \times n$ matrix with integer elements and determinant
 $\pm 1$.   Any such transformation preserves the lattice
 $\mathbb Z^n$  of integer points in $\mathbb R^n.$

As proved in Theorems \ref{theorem:volumi} and 
\ref{theorem:uniqueness}, 
up to a positive multiplicative constant, every 
rational polyhedron $P\subseteq \mathbb R^{n}$ carries a unique 
$\mathcal G_n$-invariant $d$-dimensional rational measure 
$\lambda_{d}(P), \,\,\,d=0,1,\ldots $.

To construct $\lambda_{d}$, writing  $(P,1)$ as an
abbreviation of   $\{(x,1)\in \mathbb  R^{n+1} \mid x\in P\}$,  
we  let $\Phi$ be a regular fan over the 
set    given by 
$\{\theta y \in \mathbb  R^{n+1} \mid 0\leq \theta \in \mathbb R, 
\,\,\,y \in (P,1) \}$.
Next we let  $\Delta_\Phi$ be the 
triangulation of $(P,1)$ obtained by intersecting every cone of 
$\Phi$ with the hyperplane $x_{n+1}=1.$ 
Finally, we define
 $$ 
\lambda_d(P)= \sum\left\{ \frac{1}{d!\,\,\prod_{v\in \ext(T)} 
\den(v)}\, \mid \,T \mbox{\,\,a maximal $d$-simplex of $\Delta_\Phi$} 
\right\}, $$ where $\den(v)$ denotes the least common denominator of 
the coordinates of the rational point $v\in \mathbb R^{n+1}$, and 
$\ext(T)$ is the set of vertices of $T$. 
 As a consequence of the 
celebrated solution of the weak Oda conjecture by  
Morelli and W{\l}odarczyk 
\cite{mor, wlo}, this quantity turns out not to depend on the 
chosen triangulation $\Delta_\Phi$.

As proved in Corollary \ref{corollary:lebesgue}, 
Lebesgue measure on $\mathbb R^n$  is obtainable
from $\lambda_n$ via  Carath\'eodory's construction.
However, in contrast to Lebesgue measure, for each $0\leq d<n$
and $d$-simplex  $T\in \mathscr P^{(n)},\,\,$
the rational measure  $\lambda_d(T)$  does not vanish,
and is proportional to the Hausdorff $d$-dimensional
measure of $T$,  with the constant of proportionality
only depending on the affine hull of $T$  (see
\ref{sub:pisquano}).

All necessary preliminaries on rational polyhedra (resp., cones)
and their regular triangulations  (resp., regular, or nonsingular,
 fans)  are collected in the next section.

\section{Farey blow-ups and blow-downs}
\label{section:classification}
For any rational point $x=(x_1,\ldots,x_n)\in \mathbb R^{n},
\,\,\,n=1,2,\ldots,$ we let 
$\den(x)$ denote the least common denominator of the coordinates of 
$x$.  The integer vector $$\tilde x = \den(x)(x_1,\ldots,x_n,1)\in 
\mathbb Z^{n+1}$$ is called the {\it homogeneous correspondent} of 
$x$.

For $m=0,1,\ldots$,
an  $m$-simplex $T = \conv(v_0, \ldots, v_m) \subseteq 
\mathbb R^{n}$ is said to be {\it rational} if all its vertices are 
rational.  We use the notation \begin{equation} 
\label{equation:uparrow} T^{\uparrow}= \mathbb R_{\geq 0}\,\tilde v_0 
+\cdots+ \mathbb R_{\geq 0}\, \tilde v_m \subseteq \mathbb R^{n+1} 
\end{equation} for the positive span in $\mathbb R^{n+1} $ of the 
homogeneous correspondents of the vertices of $T$.  We say that 
$T^{\uparrow}$ is the {\it (rational simplicial) cone} of $T$.  The {\it 
generators} $\tilde v_0 ,\ldots,\tilde v_m$ of $T^{\uparrow}$ are {\it 
primitive}, in the sense that each $\tilde v_i$ is minimal as a 
nonzero integer vector along its ray $\mathbb R_{\geq 0}\,\tilde v_i$.  
$T^{\uparrow}$ uniquely determines the set of its primitive 
generators, just as $T$ uniquely determines the set $ \ext(T) $ of its 
vertices.

Following \cite{ewa} we say that $T^{\uparrow}$ is {\it regular} if 
its primitive generators are part of a basis of the free abelian group 
$\mathbb Z^{n+1}$.  By definition, a rational $m$-simplex $T = 
\conv(v_0, \ldots, v_m) \subseteq \mathbb R^{n}$ is {\it regular} if 
so is $T^{\uparrow}$.

By a  {\it complex} in $\mathbb R^n$ we mean a finite set 
$\Lambda$ of compact convex polyhedra $P_i$ in $\mathbb R^n$, closed 
under taking faces, and having the further property that any two 
elements of $\Lambda$ intersect in a common face.  Unless otherwise 
specified, $\Lambda$ will be {\it simplicial}, in the sense that all 
$P_i$ are simplexes.

A simplicial complex $\Lambda$ is {\it rational} if all simplexes of 
$\Lambda$ are rational: in this case, the set $\Lambda^{\uparrow}= 
\{T^{\uparrow}\mid T\in \Lambda\}$ is a {\it simplicial fan} 
\cite{ewa, oda}.  We further say that $\Lambda$ is {\it regular} if 
the fan $\Lambda^{\uparrow}$ is regular (= nonsingular in \cite{oda}), 
meaning that every cone $T^{\uparrow}\in \Lambda^{\uparrow}$ is 
regular.

For every complex $\Lambda$, its {\it support} $|\Lambda|
\subseteq \mathbb R^n$ is the
pointset union of all polyhedra  of $\Lambda$.   Instead of
saying that $\Lambda$  is a simplicial complex with support
$P$, we briefly say that  $\Lambda$  is
a {\it triangulation} of $P$.

Given two simplicial  complexes
$\Lambda'$ and $\Lambda$ with the same support, we say that $\Lambda'$
is a {\it subdivision} of $\Lambda$ if every simplex of $\Lambda'$ is
contained in a simplex of $\Lambda$. 
For any $c \in |\Lambda|$, 
the {\it blow-up $\Lambda_{(c)}$ of
$\Lambda$ at $c$} is the subdivision of $\Lambda$ given by replacing
every simplex $C \in \Lambda$ that contains $c$ by the set of all
simplexes of the form $\conv(F\cup\{c\})$, where $F$ is any face of
$C$ that does not contain $c$ (see \cite[p.  376]{wlo}, \cite[III,
2.1]{ewa}).

 The inverse of a blow-up is called a {\it blow-down}.

 For any regular $m$-simplex $T =\conv(v_{0},\ldots,v_{m})\subseteq
 \mathbb R^{n}$, the {\it Farey mediant of $T$} is the rational point
 $v$ of $T$ whose homogeneous correspondent $\tilde v$ coincides with
 $\tilde{v}_0+\cdots+\tilde{v}_m$.  If $T$ belongs to a regular complex
 $\Delta$ and $c$ is the Farey mediant of $T$, then the
 {\it Farey blow-up}  
 $\Delta_{(c)}$ is regular.

By a {\it rational (affine)  hyperplane}  $H\subseteq \mathbb R^n$
we mean a subset of  $\mathbb R^n$ of the form
$
\{x \in \mathbb R^n \mid a\circ x +t=0\},
$
where $\circ$  denotes scalar product,  $a$ and $t$
are vectors in $ \mathbb Q^n$
(equivalently, in $\mathbb Z^n$)  and  $a\not=0.$
When $t=0$,  $H$ is called {\it homogeneous.}  

By a {\it rational affine subspace}  of $\mathbb R^n$
we mean the intersection  $A_{\mathcal F}$
of a finite set $\mathcal F$ of
rational hyperplanes in $\mathbb R^n$.
In particular,   $A_{\emptyset} = \mathbb R^n$.

The {\it affine hull}  $\,\,\aff(T)\,\,$  of  a simplex
$T$ in  $\mathbb R^n$ is the set of all
affine combinations of points of $T$.

\section{Regular triangulations  and $\mathbb Z$-homeomorphism} 

\begin{lemma} 
\label{lemma:support} 
Every rational polyhedron
$P\subseteq \mathbb R^{n}$ is the support of a regular complex. 
\end{lemma}

\begin{proof} By 
	\cite[p.36]{sta}, $P$ is the
support of some simplicial 
complex $\Lambda$.  Since $P$ is rational, $\Lambda$
can be assumed rational.  
The
set $\Lambda^{\uparrow}= \{T^{\uparrow}\mid T\in \Lambda\}$ is a
simplicial fan in $\mathbb R^{n+1}.$ The desingularization procedure
of \cite[VI, 8.5]{ewa} yields a regular subdivision $\Lambda^{*}$
of $\Lambda^{\uparrow}$.  
 Intersecting each cone of
$\Lambda^{*}$ with the hyperplane $x_{n+1}=1$ we obtain a
simplicial  complex
$\Delta$ whose support is the
set  $
(P,1)=\{(x,1)\in \mathbb R^{n+1} \mid x\in P\}. $
  For each  simplex $U\in
\Delta$ let $U'$ be the projection of $U$ onto the hyperplane
$x_{n+1}=0$, identified with
$\mathbb R^n$.  Then the regularity of $\Lambda^{*}$ ensures that
the set $\{U'\mid U\in \Delta\}$ is a regular
complex with support $P$.  \end{proof}

 The following notion is of independent interest
 \cite{buscabmun,mun-dcds}, 
 and  will find  repeated use in  this paper.

\begin{definition}
\label{definition:zhomeo}
Two rational polyhedra $P \subseteq \mathbb
R^{n}$ and $Q \subseteq \mathbb R^{m}$ are   {\it $\mathbb
Z$-hom\-e\-o\-mor\-phic,} $P\cong_{\mathbb Z} Q,$ if there is a
piecewise   linear homeomorphism 
$\eta=(\eta_1,\ldots,\eta_m)$ of $P$ onto $Q$
(each $\eta_i$  with a finite number of pieces
$l_{i1},\ldots,l_{ik(i)}$)
 such that each linear piece of
  $\eta$ and  $\eta^{-1}$ is a linear  (affine) map 
with  integer 
coefficients.     
\end{definition} 

In particular, if $m=n$  and there is $\gamma\in \mathcal G_n$
with $Q=\gamma(P)$  then 
$P\cong_{\mathbb Z} Q.$
The converse does not hold:
 the two 0-simplexes
$\{1/5\}$ and $\{2/5\}$ in $\mathbb R$
are ${\mathbb Z}$-homeomorphic
but there is no
   $\gamma\in \mathcal G_1$ 
   such that $\gamma(1/5)=2/5.$

\begin{lemma}
\label{lemma:denominator-due-claim}
Suppose
 $P \subseteq \mathbb R^{n}$ and
$P' \subseteq \mathbb R^{n'}$
are rational polyhedra and
 $\eta$ is a
  $\mathbb Z$-homeomorphism of
 $P$ onto
$P'$.

\begin{itemize}
	\item[(i)] 
	A point $z\in P$ is rational
iff so is the point $\eta(z) \in P'$.  Further, 
$
\den(y)=\den(\eta(y))\,\,\,\,\mbox{for every rational point }
y \in P.
$

	\item[(ii)] 
	There is a regular complex
	$\Lambda$ with support $P$ such that $\eta$  is linear
	(always in the affine sense) over every simplex of $\Lambda$.
	
	\item[(iii)]  For any regular complex
	$\Lambda$  with support $P$ such that 
	$\eta$  is linear
	over every simplex of $\Lambda$, the set
	$\Lambda' =  \{\eta(S)\mid s\in \Lambda\}$
	is a regular complex with support  $P'$.
	\end{itemize}
\end{lemma}

\begin{proof} 
	(i) Is an immediate consequence of
	Definition \ref{definition:zhomeo}.

(ii) { Lemma \ref{lemma:support}}  yields
a regular complex $\mathcal C_{0}$  with support $P.$
Let $\eta_{1},\ldots,\eta_{n'}$
be the components of $\eta.$
Fix  $i=1,\ldots,n'$ and let
$l_{i1},\ldots,l_{ik}$ be the
linear pieces of $\eta_{i}$.
Letting  $\sigma$ range over all  
permutations of the set   $\{1,\ldots,k\}$,
the family of  sets $P_{\sigma}=\{
x\in P\mid l_{i\sigma(1)}\leq\ldots
\leq l_{i\sigma(k)} \}$   
determines a rational  complex $\mathcal C_{i}$
with support $P$, such that
the maps $l_{ij}$  are {\it stratified}
over each polyhedron  $R$ of $\mathcal C_{i}$,
in the sense that for all $j'\not=j''$
we either have $l_{ij'}\leq l_{ij''}$ or 
$l_{ij'}\geq l_{ij''}$ on  $R$.
Since every  complex can be subdivided into
a simplicial complex without adding new vertices,
\cite[III, 2.6]{ewa}, 
we can  assume without loss of generality that all polyhedra
in $\mathcal C_{i}$ are simplexes  and
that $\mathcal C_{i}$  is a
subdivision of $\mathcal C_{0}$.
Thus $\eta_{i}$ is linear over every
simplex of $\mathcal C_{i}$ .
One now routinely constructs
a common subdivision
$\mathcal C$
of the  rational  
complexes $\mathcal C_{1},\ldots,
\mathcal C_{n'}$, such that 
every simplex of $\mathcal C$ is rational.
It follows that $\eta$  is linear
over each simplex of $\mathcal C$.
The  set $\mathcal C^{\uparrow}
 = \{T^{\uparrow}\mid T\in \mathcal C\}$
    is   a simplicial fan.
  The desingularization
   procedure  \cite[VI, 8.5]{ewa}  yields
    a regular fan $\Phi$ such that every
    cone of $\mathcal C^{\uparrow}$ is a
   union of cones of  $\Phi$.
Intersecting the cones in 
  $\Phi$ with the hyperplane $x_{n+1}=1$
we have  a   complex $\Xi$ whose support is the set
  $$(P,1)  =\{(x,1)\in \mathbb R^{n+1}\mid x\in P\}.$$
  Dropping the last
  coordinate   from the vertices of the simplexes
  of  $\Xi$ we obtain a regular 
  complex $\Lambda$ with support  $P$ such that
  $\eta$ is linear over every simplex of $\Lambda$.
   
 \medskip
  (iii)   $\Lambda'$  is a rational simplicial complex
  with support  $P'$.
Fix a  rational 
$j$-simplex  $S = \conv(v_{0},\ldots, v_{j}) 
\subseteq P$, not necessarily belonging to $\Lambda$,
such that $\eta$ is linear over $S$.
Let $S'=\eta(S)$.
The (affine) linear 
map  $\eta \colon x \in S
\mapsto y\in S'$ 
 determines the homogeneous linear map 
$  (x,1) \in (S,1) 
\mapsto (y,1) \in (S',1).
$
Let  $M_{S}$ be the $(n'+1)\times(n+1)$ integer matrix
whose bottom row has
 the form  $(0,0,\ldots,0,0,1),$
(with $n$ zeros), and
whose  $i$th row  
$\,\,\,(i=1,\ldots,n')$ is given by
the coefficients of the  (affine)  linear  
polynomial $\eta_{i}\restrict S$. 
Let   $\tilde v_{0},
\ldots\tilde v_{j}\in \mathbb Z^{n+1}$ 
be the homogeneous correspondents
of the vertices $v_{0},\ldots, v_{j}$
of $S$, and
$
S^{\uparrow} = 
\mathbb R_{\geq 0}\,\tilde v_{0}+\cdots
+\mathbb R_{\geq 0}\,\tilde v_{j}
\subseteq \mathbb R^{n+1}
$
be the positive span of
$\tilde v_{0},\ldots,\tilde v_{j}$.
 Let similarly  $S'^{\uparrow}$ be the
 positive span in $\mathbb R^{n'+1}$
 of the integer vectors
 $M_{S}\tilde v_{0},\ldots,M_{S}\tilde v_{j}$.
By construction,
$M_{S}$ sends integer points of
$S^{\uparrow}$ one-one into integer points of
$S'^{\uparrow}$. Interchanging the roles
of $S$ and $S'$ we see that
$M_{S}$ sends  integer points of
$S^{\uparrow}$ one-one {\it onto}  
integer points of
$S'^{\uparrow}$.
Blichfeldt's theorem  
  \cite[III.2]{cas},
yields the following
characterization:
\begin{align*}
 {}&  \,\, S\,\,\, \mbox{is regular}\\
 \Leftrightarrow & \mbox{ the half-open parallelepiped\,\,}
{Q}_{S} = \{\mu_{0}\tilde v_{0}+\cdots
+\mu_{j}\tilde v_{j}\mid 0\leq \mu_{0},\ldots,
\mu_{j}< 1\}\\
{}&\mbox{  contains no nonzero integer points}\\
 \Leftrightarrow&
\mbox{ the half-open parallelepiped\,\,\,}
 {Q}_{S'} \mbox{   contains no nonzero integer points}\\
  \Leftrightarrow & \,\,\, S'  \mbox{ is regular.}
\end{align*} 
In particular, if $S$ is a simplex
of $\Lambda$ then the assumed regularity of $\Lambda$ 
entails the regularity of $S$, whence
$S'$ is regular.
We conclude that  
$\Lambda'$ is a regular complex with support
$P'$.
\end{proof}

\section{Length, area, volume in $\mathcal G_n$-geometry}
For $n>0$ a fixed integer, let
  $Q \subseteq\mathbb R^n$ be a  
  (not necessarily rational) polyhedron.
For any  
triangulation $\mathcal T$ of $Q$
and $i=0,1,\ldots$
we let 
  $\mathcal T^{\max}(i)$ denote   the set of
maximal $i$-simplexes of $\mathcal T$.
The
{\it $i$-dimensional part $Q^{(i)}$ of $Q$}
is now defined by
\begin{equation}
    \label{equation:dimensional-part}
Q^{(i)}=\bigcup\{T\in \mathcal T^{\max}(i)\}.
 \end{equation}
Since 
   any two triangulations of $Q$ have a
  joint subdivision, the definition of 
$Q^{(i)}$ does not depend on the chosen
  triangulation $\mathcal T$  of $Q$.
If   $Q^{(i)}$ is nonempty, then it is  
 an $i$-dimensional polyhedron
 whose   $j$-dimensional part  $Q^{(j)}$ 
 is empty for each   $j\not= i.$
 Trivially, $Q^{(k)}=\emptyset$  for each
 integer $k>\dim(Q)$.

\smallskip
For every regular $m$-simplex  $S
=\conv(v_0,\ldots,v_m)\subseteq  \mathbb R^n$
we use the notation  
\begin{equation}
\label{equation:denominator-of-simplex}
\den(S)=\prod_{j=0}^m  \den(v_j),  
\end{equation}
and say that $\den(S)$
is the  {\it denominator} of $S$.

 \medskip 
For any   $P\in \mathscr P^{(n)}$,
regular triangulation $\Delta$ of $P$,
and $i=0,1,\ldots,$  the rational number
  $\lambda(n,i,P,\Delta)$  is defined by
  $$
\lambda(n,i,P,\Delta)= \sum_{T\in \Delta^{\max}(i)}
\frac{1}{i!\,\,\den(T)}\,,
  $$
with the proviso that
$\,\,\lambda(n,i,P,\Delta)=0\,\,$ if $\,\,\Delta^{\max}(i)=\emptyset$.
In particular, this is the case  of all $i>\dim(P).$

\medskip

\begin{proposition}
\label{proposition:triangulation-does-not-matter}
For every  $n=1,2,\ldots$,   $\, i=0,1,\ldots$, 
polyhedron
$P\in \mathscr P^{(n)}$ and regular triangulations
$\Delta$ and $ \Delta'$ of $P$, 
$
\,\lambda(n,i,P,\Delta) = \lambda(n,i,P,\Delta').
$
\end{proposition}

\begin{proof}
We first suppose that $\Delta'$ is obtained from 
$\Delta$ by a blow-up at
the Farey mediant $c$ of some $j$-simplex 
$S=\conv(v_{0},\ldots,v_{j})
\in \Delta$, $\,\,\,j=1,\ldots,n$. In symbols,
$
\Delta' = \Delta_{(c)}.
$
$S$ is the smallest simplex of $\Delta$ containing
$c$ as an element. Thus  $c\in R\in \Delta \Rightarrow \dim(R)\geq j.$ 
Let  $d=0,1,\ldots,n$.
If   for  no simplex
$T\in \Delta^{\max}(d)$ it is the case that  $c\in T$, then
$\Delta^{\max}(d)=\Delta'^{\max}(d).$ 
Otherwise, let $T=\conv(v_{0},\ldots,v_{j}, \ldots,v_{d})$  
be a   simplex  of
$\Delta^{\max}(d)$  such that  $c\in T$.
We now define the  $d$-simplexes  $S_0,\ldots,
S_j$ as follows: 
$
S_{0}=\conv(c,v_1,\ldots,v_d),
$
$
S_{j}=\conv(v_0,v_1,\ldots,v_{j-1},c,\ldots, v_d),
$
and 
$
S_{t} =\conv(v_{0},\ldots,v_{t-1}, c, v_{t+1}, \ldots,v_{j},
\ldots,v_{d})$ for each  $t=1,\ldots,j-1$.
By definition of Farey mediant, 
$\den(c) = \den(v_{0})+\cdots+\den(v_{j}).$
By {{}definition of Farey blow-up,}
the  subcomplex of $\Delta$ given by  
$T$ and its faces is replaced in $\Delta'$ by
the simplicial complex given by the $d$-simplexes
$S_{0},\ldots,S_{j}$  and their faces.
Since $T$ is regular, then so is 
$S_{u}$ for each $u=0,\ldots,j$, whence 
 $\den(S_{u})={\den(T) \cdot
\den(c)}/{\den(v_{u})}.$ 
As a consequence, 
${1}/{\den(T)}=
\sum_{u=0}^{j} {1}/{\den(S_{u})}.$
Since
$$
\sum_{T \in  \Delta^{\max}(d)} 
\frac{1}{{d!\,\den(T)}}
\,\,\,\,\,\,\,\,=\,\,\,
\sum_{U \in  \Delta'^{\max}(d)}
\frac{1}{{d!\,\den(U)}}\,,
$$
then 
$\lambda(n,d,P,\Delta)=\lambda(n,d,P,\Delta').$
Thus in case 
$
\Delta' = \Delta_{(c)} 
$ we get 
$\lambda(n,i,P,\Delta)=\lambda(n,i,P,\Delta'),$
for all $i=0,1,\ldots$.

In the general case when  $\Delta'$ is an arbitrary regular 
triangulation of $P$,
 the  {{}solution of the weak Oda conjecture} 
 \cite{mor,wlo} 
 yields 
 a sequence of regular triangulations
 $\nabla_{0}=\Delta, \,\, \nabla_{1},\ldots, 
\,\, \nabla_{s-1}, \,\,\nabla_{s}=\Delta',$ where each
 $\nabla_{k+1}$ is obtained from $\nabla_{k}$ by a   Farey
 blow-up,  or vice versa, $\nabla_{k}$ is obtained from $\nabla_{k+1}$
 by a Farey blow-up. Then the desired conclusion follows
 by induction on $s$.
\end{proof}

\medskip

In the light of the
{ foregoing proposition,}
for  each $n=1,2,\ldots,$ polyhedron
 $P\in \mathscr P^{(n)}$
and $d=0,1,\ldots,$  we can unambiguously write
\begin{equation}
\label{equation:rational-measure}
\lambda_d(P) = \lambda(n,d,P,\Delta),
\end{equation}
where  $\Delta$  is an arbitrary regular
triangulation of $P$. 
We say that $\lambda_d$  is the
$d$-dimensional  {\it rational measure}
of $P$.  Trivially,  $\lambda_d(P)=0$
for each integer  $d>\dim(P).$
If clarity requires it,  
different symbols  $\lambda_d$
and $\lambda'_d$ may be used, e.g., 
for the $d$-dimensional rational measures
defined on  $\mathscr P^{(n)}$ and
 $\mathscr P^{(n+1)}$.
 However,   no such notational distinction
 will be necessary when the   ambient space $\mathbb R^n$ 
 is  clear from the context.

 \medskip

\begin{theorem}
\label{theorem:volumi}
For  each $n=1,2,\ldots$ and 
$\,\,d=0,1,\ldots, $  the map
$\lambda_d\colon \mathscr P^{(n)}\to \mathbb R_{\geq 0}$
has the following properties,
for all $P,Q \in \mathscr P^{(n)}:$

\medskip

\begin{description}
\item[(i) \ Invariance]  If
  $P=\gamma(Q)$ for some  $\gamma\in \mathcal G_n$
then $\lambda_{d}(P)=\lambda_{d}(Q)$. 

\medskip

\item[(ii) Valuation]  $\lambda_{d}(\emptyset)=0,\,\,$
$\lambda_{d}(P)=\lambda_d(P^{(d)})$,  and the restriction of 
 $\lambda_{d}$ to the set of all 
 rational polyhedra  $P,Q$  in $ \mathbb R^{n}$ 
 having dimension at most $d$ is a {\rm valuation}: in other words,
\begin{equation}
\label{equation:valuation}
		\lambda_{d}(P)+\lambda_{d}(Q)=
			\lambda_{d}(P\cup Q)+
			\lambda_{d}(P\cap Q).
\end{equation}

 \medskip

\item[(iii) Conservativity]  For any $P\in \mathscr P^{(n)}$
let $(P,0)=\{(x,0)\in \mathbb R^{n+1}\mid x\in P \}.$
Then $\lambda_d(P)=\lambda_d(P,0)$. 

\medskip

\item[(iv) Pyramid]   For
$k=1,\ldots,n$, if
  $\,\,\conv(v_0,\ldots,v_k)$
is a regular $k$-simplex in $\mathbb R^n$
with  $v_0\in \mathbb Z^n$ then
\begin{equation}
\label{equation:pyramid}
\lambda_k(\conv(v_0,\ldots,v_k))=
 {\lambda_{k-1}(\conv(v_1,\ldots,v_k))}/{k}.
\end{equation}
  
\medskip
  
\item[(v) Normalization]
Let 
$j=1,\ldots,n$.  Suppose  the set 
  $B=\{w_1,\ldots,w_j\} \subseteq  \mathbb Z^n$
is part of a basis of the free abelian 
  group $\mathbb Z^n$.
  Let the closed   
parallelepiped $P_B\subseteq \mathbb R^n$
  be defined by
\begin{equation}
\label{equation:parallelepiped}
  P_B =\left\{x\in \mathbb R^n\mid x=
  \sum_{i=1}^j \gamma_iw_i,\,\,\, 0\leq \gamma_i\leq 1\right \}.
\end{equation}
  Then  $\lambda_j(P_B)=1.$
  
  \medskip
  
  \item[(vi) Proportionality]  Let 
$A$ be an $m$-dimensional rational affine subspace of $\mathbb
 R^n$  for some $m=0,\ldots,n$.  
 Then there is a constant $\kappa_A>0$,
 only depending on $A$,   such that 
  $\,\,\,\lambda_m(Q) = \kappa_A\cdot \mathcal H^m(Q)$
  for every rational 
  $m$-simplex
   $Q\subseteq A.$
Here as usual,  $\mathcal H^m$ denotes $m$-dimensional
 Hausdorff measure. 
\end{description}
\end{theorem}

In Theorem \ref{theorem:uniqueness}  below we will prove that
conditions (i)-(vi)  uniquely characterize the maps $\lambda_d
\colon \mathscr P^{(n)}\to \mathbb R_{\geq 0}.$
We first prove that  conditions
(i)-(vi) hold.

\section{Proof of Theorem \ref{theorem:volumi}(i)-(v)}
\subsection{Invariance}
We will actually prove the 
stronger result that 
$\lambda_d$ is invariant under $\mathbb Z$-homeomorphisms:
in other words, whenever 
$P'\subseteq \mathbb R^{n'}$ is a rational polyhedron and
$P\cong_{\mathbb Z} P'$ then $\lambda_{d}(P)=\lambda_{d}(P')$ for all
$d=0,1,\ldots$.
Let $\iota$
be a  $\mathbb Z$-homeomorphism
of $P$ onto $P'$.
Let  $\Delta$ be 
a regular complex with support 
$P$ such that  
$\iota$ is  (affine) linear
over every simplex of $\Delta$.
The existence of $\Delta$  is ensured by
{ Lemma \ref{lemma:denominator-due-claim}(ii).}
Let $\Delta'= \{\iota(T)\mid T
\in \Delta\}$.
By {Lemma \ref{lemma:denominator-due-claim}(i)-(iii),}
$\Delta'$  is a regular
complex with support $P'$, and
$
\den(\iota(z)) = \den(z)$
for every rational point
$z\in P.$
It follows that
$
\lambda(n,d,P,\Delta)=
\lambda(n',d,P',\Delta').
$
The desired conclusion now follows
from (\ref{equation:rational-measure})
as a consequence of
{Proposition \ref{proposition:triangulation-does-not-matter}.}

 
\subsection{Valuation}
The identities
$\lambda_{d}(\emptyset)=0,\,\,$  and 
$\lambda_{d}(P)=\lambda_d(P^{(d)})$ 
immediately follow by 
{ definition  of rational measure}.
To prove (\ref{equation:valuation}), 
 we first observe that both $P\cup Q$ and
$P\cap Q$ 
are rational polyhedra in $\mathbb R^{n}$
		whose dimension is at most $d$.
As an application of
{Lemma \ref{lemma:support},} let the
regular complexes
$\Delta, \Phi, \Psi, \Omega$  
have the following properties:  
$$
|\Delta|=P\cap Q, \,\,\,|\Phi|=P,
\,\,\,|\Psi|=Q, \,\,\,|\Omega|=P\cup Q.
$$
Using   the   extension argument in
\cite[VI. 9.3]{ewa}  we can 
 assume 
$\Delta=\Phi\cap \Psi$
and $\Omega  =  \Phi\cup \Psi$,
without loss of generality. 
For every $X\subseteq \mathbb R^{n}$ we
let as usual  $\cl(X)$
denote the closure of $X$ in $\mathbb R^{n}$.
By { Proposition \ref{proposition:triangulation-does-not-matter}}
 we have:

\[
\begin{array}{lll}
& &\lambda_{d}(P)+\lambda_{d}(Q) = 
\lambda(n,d,P,\Phi)+\lambda(n,d,Q,\Psi)\\[0.4cm]
&=& 
\frac{1}{d \,!}\left[\sum_{T \in  \Phi^{\max}(d)}{\den(T)}^{-1}
\,+\,\,\,
\sum_{T \in  \Psi^{\max}(d)} {\den(T)}^{-1}\right]  \\[0.4cm]

&=& 
\frac{1}{d \,!} 
\left[
\sum_{\cl(P\setminus Q) \supseteq T \in  \Phi^{\max}(d)} {\den(T)}^{-1} 
+ \sum_{\cl(Q\setminus P) \supseteq T \in
  \Psi^{\max}(d) } {\den(T)}^{-1} 
  \right] 
 \\[0.4cm]

&  &   +
\frac{1}{d \,!} 
\left[ \sum_{{P\cap Q\supseteq  T
 \in  \Phi^{\max}(d)}} {\den(T)}^{-1} 
 +
\sum_{ {P\cap Q \supseteq  T \in  
\Psi^{\max}(d)} } {\den(T)}^{-1}
\right]\\[0.4cm]

&=& 
\frac{1}{d \,!}
\left[
\sum_{ {\cl(P\setminus Q) \supseteq T \in  
\Phi^{\max}(d)}} {\den(T)}^{-1}
 +
\sum_{{\cl(Q\setminus P)
 \supseteq  T \in  
\Psi^{\max}(d)}} {\den(T)}^{-1}
\right] \\[0.4cm]

& & +  
\frac{2}{d \,!}\sum_{{ T \in  \Delta^{\max}(d)}}  {\den(T)}^{-1} 

 \end{array}
 \]
 
 \[
\begin{array}{lll}

&=&  \frac{1}{d \,!}
\left[\sum_{ {\cl(P\setminus Q) \supseteq T \in  
\Omega^{\max}(d)}} {\den(T)}^{-1}
+ \sum_{{\cl(Q\setminus P)  \supseteq  T \in  
\Omega^{\max}(d)}} {\den(T)}^{-1}
\right]  \\[0.4cm]
& &  +  
\frac{2}{d \,!}\sum_{ { T \in  \Delta^{\max}(d)}}  {\den(T)}^{-1} \\[0.4cm]

&=&
\frac{1}{d \,!}
\left[
\sum_{{\cl(P\setminus Q) \supseteq T \in  
\Omega^{\max}(d)}}  {\den(T)}^{-1}
+ \sum_{{\cl(Q\setminus P)  \supseteq  T \in  
\Omega^{\max}(d)}}  {\den(T)}^{-1}
\right] \\[0.4cm]
& &  + 
\frac{1}{d \,!}
\left[ 
\sum_{{P\cap Q \supseteq  T \in  
\Omega^{\max}(d)}}  {\den(T)}^{-1}
 +
\sum_{{ T \in  \Delta^{\max}(d)}}  {\den(T)}^{-1}
\right] \\[0.4cm]

&=& \lambda(n,d,P\cup Q,\Omega) +
\lambda(n,d,P\cap Q, \Delta) =
  \lambda_{d}(P\cup Q)+\lambda_{d}(P\cap Q).
\end{array}
\]

\smallskip
\subsection{Conservativity and Pyramid}
 
Properties  (iii)-(iv)
are immediate consequences of the definition
of $\lambda_d.$

\subsection{Normalization} 
\label{sub:normalization}
To prove Property (v),  
let  $\Pi$ be the set of permutations
 of the set $\{1,2,\ldots,j\}$. For every permutation
 $\pi\in\Pi$ we let  $T_\pi$ be the convex hull of
 the set of points  
 $$0,\,\,w_{\pi(1)},\,\,w_{\pi(1)}+w_{\pi(2)},\,\, w_{\pi(1)}+w_{\pi(2)}+w_{\pi(3)},
 \ldots,w_{\pi(1)}+w_{\pi(2)}+\cdots+
 w_{\pi(j)}. $$
Arguing as in 
  \cite[3.4]{sem},  it follows that
  the $j$-simplexes $T_\pi$ are the maximal
  elements of a  triangulation $\Sigma$
 of  $P_B$, called the
 {\it standard triangulation}  $\Sigma$. Each simplex
 $T_\pi$  is regular and 
 and has unit denominator. There are $j!$ such simplexes.
 By definition,  the rational 
 $j$-dimensional measure
 of  $T_\pi$   is  equal to  $1/{j!}\,\,$.
 A final application of   
{\ref{proposition:triangulation-does-not-matter} } yields   
 $\lambda_j(P_B)=1.$

\section{From $\lambda_n$ to Lebesgue measure on
$\mathbb R^n$ via Carath\'eodory's method}
In what follows, 
   $\mathcal L^n$ will
 denote  Lebesgue measure on $\mathbb R^n$.

\begin{proposition} 
     \label{proposition:lebesgue-generalizzato} 
    For any  $n=1,2,\ldots$ and polyhedron  
$Q\in  \mathscr P^{(n)}$, 
$\lambda_n(Q)= \mathcal L^n(Q).$
\end{proposition}

\begin{proof} If $\dim(Q)<n$ then $\mathcal L^n(Q)=\lambda_n(Q)=0$.  
If $\dim(Q)=n$, since
$\lambda_n(Q)=\lambda_n(Q^{(n)})$ and $\mathcal L^n(Q)=
\mathcal L^n(Q^{(n)})$,
 without loss of generality we may assume $Q=Q^{(n)}$.
 Let $\nabla$ be a regular triangulation of $Q$ as given by { Lemma 
\ref{lemma:support}.} 
Since, as we have seen,  $\lambda_n$ is a valuation on
$\mathscr P^{(n)}$ and $\mathcal L^n(Q)=
\sum_{S\in \nabla^{\max}(n)} \mathcal L^n(S)$, it 
is enough to prove 
\begin{equation}
\label{equation:desideratum}
\lambda_n(S)= \mathcal L^n(S)\,\,
\mbox{ for every } n\mbox{-simplex }
S=\conv(w_{0}, \ldots, w_{n})\in \nabla.
\end{equation}
To this purpose, 
let $T \subseteq {\mathbb R}^{n+1}$ be the $(n+1)$-simplex with 
vertices $0$, $(w_{0},1)$, $\ldots$, $(w_{n},1)$.  Then 
$$
 \mathcal L^{n+1}(T) = {\mathcal L^n({S})}/({n+1}).  
 $$
  This is the classical 
formula for the volume of the $(n+1)$-dimensional pyramid with base 
$S$ and height 1.
Next we observe that
$T$ is contained in the closed $(n+1)$-dimensional parallelepiped 
$E = 
\{\alpha_{0}(w_{0},1)+ \cdots+\alpha_{n}(w_{n},1) \in {\mathbb 
R}^{n+1} \mid \alpha_{0},\ldots,\alpha_{n} \in [0,1]\}.
$
Further, $E\subseteq { U} = \{\alpha_{0}{\tilde w_{0}}+ 
\cdots+\alpha_{n}{\tilde w_{n}} \in {\mathbb R}^{n+1} \mid 
\alpha_{0},\ldots,\alpha_{n} \in [0,1 ]\}.
$
 Since $S$ is regular, a 
classical argument in the geometry of numbers, (\cite{cas} or 
\cite[Proof of VI, 8.5]{ewa}) yields $\mathcal L^{n+1}({ U})=1$. 
 For all $i=0,\ldots,n$ let $d_i=\den(w_i)$.  
Since ${\tilde w_{0}} = d_{0} (w_{0},1),\ldots, {\tilde w_{n}} = d_{n} 
(w_{n},1)$, then $$\mathcal L^{n+1}(E) = (d_{0} \cdots d_{n})^{-1}.$$
The construction of \cite[3.4]{sem}  now  yields
 a triangulation of $E$ 
consisting of $(n+1)$-simplexes $T_1,\ldots,T_{(n+1)!}$ and their faces, 
in such a way that 
 $$ 
\mathcal L^{n+1}(T_i) = \frac{\mathcal L^{n+1}(E)}{(n+1)!}
\,\,\,\mbox{ for each 
$i=1,\ldots,(n+1)!$}  \,\,\, 
$$  
Each  $T_i$ is a regular simplex.
One easily gets  a  linear  
(affine) isometry  of $T_i$ onto $T$.  Therefore,   
$$ \mathcal L^{n+1}(T) = 
\frac{\mathcal L^{n+1}(E)}{(n+1)!}\,.  
$$
 Summing up, $\mathcal  L^n({S}) = 
 \mathcal L^{n+1}(E)/n!= ({n!  \, d_{0} \cdots d_{n}})^{-1} 
= \lambda_n(S)$, 
and (\ref{equation:desideratum}) is proved.  
\end{proof}

\begin{corollary}
\label{corollary:lebesgue} 
Fix $n=1,2,\ldots$ and let
$\mathscr K^{(n)}$ denote the family
of compact subsets of $ \mathbb R^n$. 
For any Borel set  $E\subseteq \mathbb R^n$ let us define
$$
\bar\lambda_n(E)=
\sup_{E\supseteq K\in \mathscr K^{(n)}}\,\,
\inf_{K\subseteq P \in  \mathscr P^{(n)}}
\lambda_n(P).
$$
Then $\bar\lambda_n(E)=\mathcal L^n(E).$
\end{corollary}

\begin{proof}  We first {\it claim} that
every   $K\in \mathscr K^{(n)}$
coincides with the intersection of all rational polyhedra of 
$\mathscr P^{(n)}$ containing it.

As a matter of fact, 
for  any  $P,Q\in  \mathscr P^{(n)}$
both $P\cup Q$ and $P\cap Q$ are members of $\mathscr P^{(n)}$.
Moreover, there  exists a rational triangulation
 $\mathcal T$  of $P\cup Q$ such that the set
 $\{T\in \mathcal T\mid T\subseteq P\cap Q\}$
 is a triangulation of   
 $P\cap Q$. Thus the set
  $\{T\in \mathcal T\mid T\subseteq \cl(P\setminus Q)\}$
  is a triangulation of  the set  $\cl(P\setminus Q)\subseteq 
  \mathbb R^n$, which shows that 
   $\cl(P\setminus Q) $ is a rational polyhedron. 
 For every $x\in \mathbb R^n\setminus K$ there is
a rational $n$-simplex  $T$ containing $x$ in its interior
and such that $T\cap K=\emptyset.$
Since $K$ is contained in some rational polyhedron, 
our claim is settled.

 Now let $P_0\supseteq P_1\supseteq \cdots$
be a sequence of rational polyhedra such that
$\bigcap_i P_i=K$ and for every   $R \in \mathscr P^{(n)}$ with
$K\subseteq R$ there is
$j=0,1,\ldots$ such that $P_j\subseteq R.$
The existence of this sequence follows
from our claim together with the observation that there are 
 only countably many   rational polyhedra.
By   {Proposition \ref{proposition:lebesgue-generalizzato}}, 
 $\lambda_n(P_0)=\mathcal L^n(P_0)\geq\mathcal L^n(P_1)
 =\lambda_n(P_1)\geq\lambda_n(P_2)\geq\cdots,$  whence by 
  construction, 
$ \lim_{i\to \infty} \lambda_n(P_i)=
\inf\{\lambda_n(R)\mid R 
\supseteq K,\,\, R \in \mathscr P^{(n)}\}=\bar\lambda_n(K).$
{}Combining
{Proposition \ref{proposition:lebesgue-generalizzato}}
with   the countable monotonicity property of $\mathcal L^n$,
 we get
 $\mathcal L^n(K)=\lim_{i\to \infty} \mathcal L^n(P_i)=
 \lim_{i\to \infty} \lambda_n(P_i)=\bar\lambda_n(K).$  
 
 Having thus proved that  $\bar\lambda_n$
 agrees with $\mathcal L^n$ on all compact subsets of
 $\mathbb R^n$, the desired conclusion
   follows from the regularity  properties of Lebesgue measure.
 \end{proof}

\begin{remark}
Following   \cite[115C]{fre}  we 
 now routinely extend  $\bar\lambda_n$
to an outer measure  $\lambda^{*}_n\colon 
 \mbox{powerset}({\mathbb R^n})
\to [0,\infty]$ which, by { Corollary \ref{corollary:lebesgue}
and \cite[115D]{fre}}
coincides with Lebesgue outer measure on $\mathbb R^n$.   
As proved in
   \cite[115E]{fre},  by applying to $\lambda^{*}_n$   
Carath\'eodory's construction \cite[113]{fre}
we finally obtain
Lebesgue measure on $\mathbb R^n$.
\end{remark}

\section{Proof of Theorem \ref{theorem:volumi}(vi)}
\subsection{Basic material on Hausdorff measure}
In the following proposition we collect
a number of well known consequences of the
isodiametric inequality
(see \cite[2.10.33]{fed}), 
and of the invariance of Hausdorff 
$d$-dimensional measure
under isometries:

\begin{proposition}
    \label{proposition:natural-measure}
For each $0<n\in \mathbb Z$
we have:

 \medskip

\begin{itemize}
\item[(i)] 
If $\,T=\conv(x_0,\ldots,x_n)$ is an $n$-simplex in
$\mathbb R^n$, letting
  $M$ be the $n\times n$ matrix whose $i$th
row is given by the vector $x_i-x_0,\,\,\,
(i=1,\ldots,n),$  it follows that 
 $\mathcal H^n(T) = |\det(M)|/n! = \mathcal L^n(T).$

\medskip
\item[(ii)]  If $S$ is an $m$-simplex in
$\mathbb R^n$
with $0< m<n$,  and we map  $S$
onto a   copy $S'$ by means
of an
isometry $\iota$ sending the
affine hull of  $S$
onto the linear 
subspace $\mathbb R^m$ of $\mathbb R^{n}$
spanned by the first $m$ standard basis vectors
of $\mathbb R^n$,
then 
$\mathcal H^m(S)=\mathcal L^m(S')$.
If $\dim(S)=0$ then
$\mathcal H^0(S)=1=$ number of elements of
the singleton $S$.

\medskip
\item[(iii)]  
 If $Q$ is  a  nonempty polyhedron  in   $\mathbb R^n$
 and  $Q=Q^{(d)}, \,\,\,d=0,1,\ldots$,  then 
 letting  $\mathcal T$ be an arbitrary triangulation 
 of  $Q$,
with its $d$-simplexes 
$T_1,\ldots,T_k$,   we have
$\mathcal H^{d}(Q)= 
\sum_{j=1}^k \mathcal H^{d}(T_j).$
If 
$Q=\emptyset$ then  $\mathcal H^{k}(Q)=0$ for all
$k=0,1,\ldots$.      

\medskip
\item[(iv)] 
   Given integers  $0\leq m < n$, 
suppose 
$T=\conv(v_{0},\ldots,v_m)$ and
$T'=\conv(v'_{0},\ldots,v'_m)$
are  $m$-simplexes  in $ \mathbb R^n$ 
with 
$\,\aff(T) = \aff(T')$.  For $v$  
an arbitrary point lying in
$\mathbb R^{n}\setminus \aff(T)$,   let $ U=\conv(T,v)$ and $
U'=\conv(T',v).  $ 
Then 
$ {\mathcal H^{m+1}(U')}/{\mathcal H^{m+1}(U)} = 
{\mathcal H^m(T')}/{\mathcal H^m(T)}.  $

\medskip
\item[(v)]  
More generally,  suppose  the points  
$v_{m+1},\ldots,v_n \in \mathbb R^n$
have the property that 
$W=\conv(v_{0},\ldots,v_m,
v_{m+1},\ldots,v_n)$  is
an $n$-simplex.  Then also 
$W'=\conv(v'_{0},\ldots,v'_m,
v_{m+1},\ldots,v_n)$ is an
$n$-simplex, and  we have the identity
${{\mathcal H^n(W')}}/{
{\mathcal H^n(W)}} = 
{{\mathcal H^m(T')}}/{{\mathcal H^m(T)}}.  
$

\end{itemize}
\end{proposition}


\medskip 
\subsection{End of the proof of Theorem \ref{theorem:volumi}(vi)}
\label{sub:pisquano}
There remains to be proved that  $\lambda_d$ has
the Proportionality property (vi).
By { Lemma \ref{lemma:support}}, 
 $Q$ has a
regular triangulation.  
Since  $\lambda_m$ is a valuation,
recalling {Proposition \ref{proposition:natural-measure}(iii)}
 it
 suffices to argue in case $Q$ is a
regular $m$-simplex.  If $m=n$  the result follows 
from  {Proposition 
\ref{proposition:lebesgue-generalizzato}} since, 
by  {Proposition 
\ref{proposition:natural-measure}(i),} 
$\mathcal H^n(Q)=\mathcal L^n(Q)$.
In this
case $\kappa_A=1.$
Next suppose $0\leq m<n.$ It suffices to prove that for any two regular
$m$-simplexes $T =\conv(v_0,\ldots,v_m)$ and
$T'=\conv(v'_0,\ldots,v'_m)$ lying in  $A,$ 
$${\lambda_m(T)}/{\lambda_m(T')}= {\mathcal H^m(T)}
/{\mathcal H^m(T')}.  $$ 
To this purpose let $U=
\conv(v_0,\ldots,v_m,v_{m+1},\ldots,v_n)$ be a regular $n$-simplex
in $\mathbb R^n$ having $T$ as a face.

\medskip
\noindent
{\it Claim.}  
The simplex
$U'=\conv(v'_0,\ldots,v'_m,v_{m+1},\ldots,v_n)$
is regular. 

\smallskip
 
As a matter of fact, the regularity of $T$  means that the set
$\{\tilde v_0,\ldots,\tilde v_m\}$ is a basis of the free abelian 
group $G = \mathbb Z^{n+1} \cap (\mathbb R{\tilde v_0}+\cdots+\mathbb 
R{\tilde v_m})$ of integer points in the $(m+1)$-dimensional linear 
space spanned by $\tilde v_0,\ldots,\tilde v_m$ in $ \mathbb R^{n+1}$.  
Since $\aff(T')=A=\aff(T)$ and $T'$ is regular, also $\tilde 
v'_0,\ldots,\tilde v'_m$ constitute a basis of $G$.  Upon writing each 
$\tilde v_i$ and $\tilde v'_j$ as a column vector, let $M$ be the 
$(n+1)\times(m+1)$ matrix whose $i$th row coincides with $\tilde v_i$.  
Let similarly $M'$ be the $(n+1)\times(m+1)$ matrix whose $j$th row 
equals $\tilde v'_j.$ Let the $(m+1)\times(m+1)$ integer matrix $Z$ be 
defined by $MZ=M'$.  The $(m+1)\times(m+1)$ integer matrix $V$ defined 
by $M'V=M$ coincides with $Z^{-1}$, whence 
$|\det(Z)|=|\det(Z^{-1})|=1.$ Let the matrix $N$ be defined by $$ N 
=\left(\begin{tabular}{c|c} $Z$ & $0$ \\ \hline $0$ & ${\rm I}_{n-m}$ 
\end{tabular}\right) $$ where ${\rm I}_{n-m}$ denotes the $(n-m)\times 
(n-m)$ identity matrix.  $N$ is a unimodular integer 
$(n+1)\times(n+1)$ matrix.  Let $W$ (resp., $W'$) be the 
$(n+1)\times(n+1)$ integer matrix whose first $m+1$ columns are those 
of $M$ (resp., those of $M'$), and whose last $n-m$ columns are given 
by the column vectors $\tilde v_{m+1},\ldots,\tilde v_{n}.$ {}From 
$WN=W', $ it follows that the vectors $\tilde v'_0,\ldots,\tilde v'_m, 
\tilde v_{m+1},\ldots,\tilde v_n$ constitute a basis of the free 
abelian group $\mathbb Z^{n+1}$.  Therefore, 
$\,\,\,\conv(v'_0,\ldots,v'_m, v_{m+1},\ldots, v_{n})\,\,\,$ is a 
regular $n$-simplex in $\mathbb R^n$, and our claim is settled.

\medskip
Let now
$d_i=\den(v_i),\,\,\, (i=0,\ldots,n)$ and $d'_j=\den(v'_j),\,\,\,
(j=0,\ldots,m).$  Since both simplexes $U$ and $U'$ are regular
we can write the identities

\medskip
$$
\frac{\lambda_m(T)}{\lambda_m(T')} 
=\frac{(m!\,\,d_0\cdots d_m)^{-1}}{(m!\,\,d'_0\cdots d'_m)^{-1}}
=\frac{(n!\,\,d_0\cdots d_md_{m+1}\cdots d_n)^{-1}}
{(n!\,\,d'_0\cdots d'_md_{m+1}\cdots d_n)^{-1}}
=\frac{\lambda_n(U)}{\lambda_n(U')}.
$$

\medskip
\noindent By { Propositions
\ref{proposition:lebesgue-generalizzato}
and     
\ref{proposition:natural-measure}(ii)-(v)} we obtain

\smallskip
$$
\frac{\lambda_n(U)}{\lambda_n(U')}
=\frac{{\mathcal L}^n(U)}{{\mathcal L}^n(U')}
=\frac{{\mathcal H}^n(U)}{{\mathcal H}^n(U')}
=\frac{{\mathcal H}^m(T)}{{\mathcal H}^m(T')}\, ,
$$

\smallskip
\noindent
as required to prove  (vi).

 \bigskip

 The proof of Theorem \ref{theorem:volumi}  is now complete. 
 \hfill $\Box$

\section{Uniqueness}
For every nonempty rational affine subspace $F$ of $\mathbb R^n$ let the 
integer $d_F \geq 1$ be defined by 
\begin{equation} 
\label{equation:d} 
d_F = \min\{q \in \mathbb Z \mid q=\den(r)
 \mbox{\,\,\,for some 
rational point} \,\,\,r\in F\}.  
\end{equation}

\begin{lemma}
\label{lemma:wlog-isobunch}  Fix $n=1,2,\ldots$ and $e=0,\ldots,n$.
Let $F$ be a rational 
$e$-dimensional affine subspace of
$\mathbb R^n$ and $\,\,d=d_F$.
\begin{itemize}
 \item[(i)] There are rational points
 $\,\,v_0,\ldots,v_e\in F$, all with denominator $d$,
 such that $\conv(v_0,\ldots,v_e)$ is a regular
 $e$-simplex. 
 
  \item[(ii)]
   For any rational point $y\in F$
 there is an integer
  $\,\,k =1,2,\ldots$  such that $\den(y)=k d.$
 \end{itemize}
\end{lemma}

\begin{proof}
(i)  For some regular $e$-simplex  
$S_0= \conv(u_0,\ldots,u_e)$  we can write  
$ F=\aff(u_0,\ldots,u_e)$.
The regularity of $S_0$ means that the set
 $B_0=\{\tilde u_0,\ldots,\tilde u_e\}$  can be extended to a basis
 of the free abelian group $\mathbb Z^{n+1}$, whence $B_0$ is
 a basis of the lattice  $\mathbb Z^{n+1}\cap F^*, $
 where
 $F^*=\mathbb R\tilde u_0+\cdots +\mathbb R\tilde u_e$
is  the  linear subspace of
 $\mathbb R^{n+1}$  generated by 
 $ \tilde u_0,\ldots,\tilde u_e.$
 
It is impossible that the {\it heights}  (=last coordinates)
of $\tilde u_0,\ldots,\tilde u_e$ are all equal to an integer
$h>d$:  for otherwise, 
any  primitive vector $\tilde r$ in $F^*$  of height $d$,
for $r$ as in  (\ref{equation:d}),  could not arise
as a linear combination of the $\tilde u_i$ with
integer coefficients---and $B_0$ would not be
a basis of $\mathbb Z^{n+1}\cap F^*$.

If  the heights of  
$\tilde u_0,\ldots,\tilde u_e$  are all equal to $d$
we have nothing to prove.
Otherwise,  we will construct a finite sequence $B_0, B_1,\ldots$
 of bases of 
   $\mathbb Z^{n+1}\cap F^*$, 
 and  finally  obtain a  basis 
 $\{\tilde v_0,\ldots,\tilde v_e\}$
 having the property that the
 height of each  $\tilde v_i$  is equal to $d$.

The first step is as follows:   
Choose a vector $\tilde u_i\in B_0$  of top height, a vector 
$\tilde u_j\in B_0$  of smaller height, and replace
$\tilde u_i$  by $\tilde u_{i}-\tilde u_j.$  We get a new basis
$B_1$  of  $\mathbb Z^{n+1}\cap F^*$, and
a new regular $e$-simplex  $S_1$  in $F$.
Specifically, letting the rational point $w\in F$
be defined by $\tilde w = \tilde u_{i}-\tilde u_j,$
the vertices of $S_1$ are $u_0,\ldots,u_{i-1},
w, u_{i+1},\ldots,u_e$. Observe that the sum of the
heights of the elements of $B_1$
is strictly smaller than the sum
of the heights of the elements of $B_0$.

Proceeding inductively, and replacing a top
vector $\tilde u$ of the basis $B_t$ by a vector $\tilde u-\tilde v$
with $\tilde v\in B_t$ of smaller height than $\tilde u$, we obtain
a new basis $B_{t+1}$ such that the 
sum of the heights of the elements  of $B_{t+1}$ is strictly smaller
than the sum of the heights of the elements  of $B_{t}$. 
We also get a new regular $e$-simplex  $S_{t+1}$ lying in $F$.
The process must terminate with a basis
$\{\tilde v_0,\ldots,\tilde v_e\}$ of
$\mathbb Z^{n+1}\cap F^*$ where all
$\tilde v_i$ have the same height---which by
our initial discussion must be equal to $d$.
By definition, $\conv(v_0,\ldots, v_e)$ is the desired regular
$e$-simplex in $F$.

\smallskip
(ii) now trivially follows from (i): 
for,  the   regularity of
$\conv(v_0,\ldots,v_e)$  is to the effect that
 the primitive vector
$\tilde y\in \mathbb Z^{n+1}$ is a linear combination of the
$\tilde v_i$ with integer coefficients.
\end{proof}

\begin{theorem}
\label{theorem:uniqueness}
Properties (i)-(vi) in Theorem
\ref{theorem:volumi} uniquely characterize the rational
measures
$\lambda_0,\ldots,\lambda_n,$  for each
$n=1,2,\ldots$,  among all
maps from
$ \mathscr P^{(n)}$  to  $\mathbb R_{\geq 0}.$ 
\end{theorem}

\begin{proof}  Suppose
for each  $n=1,2,\ldots$, the maps 
$\mu_0,\ldots,\mu_n\colon 
 \mathscr P^{(n)}\to \mathbb R_{\geq 0}$,
 as well as the maps
 $\mu'_0,\ldots,\mu'_{n+1}\colon 
 \mathscr P^{(n+1)}\to \mathbb R_{\geq 0}$ 
have all properties  (i)-(vi).
 Since by {Lemma \ref{lemma:support}}
every rational polyhedron has a regular triangulation,
and  each $\mu_j$ and $\lambda_j$  is a { valuation},
it suffices to show that
$\mu_m(S)=\lambda_m(S)$  for all 
 $m=0,\ldots,n,$  and regular
$m$-simplex $S$ in $\mathbb R^n.$
Let $F=\aff(S)$  be the affine hull of $S$
in $\mathbb R^n.$  Let $d=d_F$ be the  smallest
denominator of a rational point of $F$  as in
{ (\ref{equation:d})  above.}
Let  us identify $\mathbb R^n$
with the hyperplane  $x_{n+1}=0$ of $\mathbb R^{n+1}$. 
Let $T=\conv(v_0,\ldots,v_m)\subseteq F$ 
be a regular $m$-simplex  such that
$\den(v_0)=\cdots=\den(v_m)=d.$  The existence of
$T$ is ensured by { Lemma \ref{lemma:wlog-isobunch}.} Let
$T'=\{(x,1)\in \mathbb R^{n+1}\mid x\in T \}$.
There is  $\alpha\in \mathcal G_{n+1}$ such that
$\alpha(T,0)=T'.$
{} From the
 {Invariance and Conservativity} properties of
all $\mu_i$  and $\mu'_j$  we obtain
\begin{equation}
\label{equation:lorena}  
\mu'_m(T')= \mu'_m(T,0)=\mu_m(T).
\end{equation}
The regularity of $T$ means  that the set
$B=\{\tilde v_0,\ldots,\tilde v_m\}$
is part of a basis of the free abelian group
$\mathbb Z^{n+1}$. 
As  in 
(\ref{equation:parallelepiped})  above, let the
 closed
parallelepiped $P_B$  be defined by
$$
P_B =\left\{x\in \mathbb R^{n+1}\mid x=
  \sum_{i=0}^m \gamma_i\tilde v_i,\,\,\, 0\leq \gamma_i\leq 1\right \}.$$
{}From the {Normalization} property
we get    $\mu'_{m+1}(P_B)=1.$
Arguing as in   
 \ref{sub:normalization} above,  we obtain a
 triangulation  $\Delta$ of $P_B$   consisting of 
 $(m+1)$-simplexes  $T_1,\ldots,T_{(m+1)!}$ and their faces.
Each  $T_i$ is
regular
and has denominator 1.  A direct verification
shows that for any two such simplexes
$T_i$  and $T_j$  there is $\gamma\in \mathcal G_{n+1}$
such that $T_i=\gamma(T_j).$ From the
{Valuation and Invariance}  properties of
$\mu'_{m+1}$
 it follows that
 $$
 \mu'_{m+1}(T_j)=  \frac{\mu'_{m+1}(P_B)}{(m+1)!} = \frac{1}{(m+1)!}
 \quad \mbox{for all }  j=1,\ldots,(m+1)! 
 $$
 Let  $D\subseteq \mathbb R^{n+1}$  be the $(m+1)$-simplex
 with vertices $0, \tilde v_0,\ldots,\tilde v_m.$
 It is easily seen that $D$ is regular and  $\den(D)=1.$ 
 Thus an easy exercise yields 
  an  $\eta\in  \mathcal G_{n+1}$ such that
 $\eta(T_1)=D.$  One more application of
  the {Invariance} property
of $ \mu'_{m+1}$ yields 
 $$
 \mu'_{m+1}(D)= \frac{1}{(m+1)!}\,.
 $$
Since the $(m+1)$-simplex  $D'$  with vertices
 $0,(v_0,1),\ldots,(v_m,1)$ has the same affine hull as
 $D$, by the assumed {Proportionality} 
 property of $ \mu'_{m+1}$  we have
 $$
 \mu'_{m+1}(D') = \frac{1}{(m+1)!\,d^{m+1}}\,.
 $$
On the other hand,  the {Pyramid} property 
is to the effect that
 $$
  \mu'_{m+1}(D') = \frac{ \mu'_{m}(T')}{m+1},
 $$
 whence
 $$
 \mu'_{m}(T')=\frac{1}{m!\,d^{m+1}}\, .
 $$ 
 Recalling  (\ref{equation:lorena})
 we get
 $$
  \mu_{m}(T)=\frac{1}{m!\,d^{m+1}} =\lambda_m(T),
 $$
 because $T$ is regular and all the denominators 
 of its vertices are
 equal to $d$.
Since $S$ and $T$ have the
same affine hull, a
 final application of the {Proportionality property}
 of $\mu_{m}$  and $\lambda_{m}$
  yields
 $$
 \frac{\mu_{m}(S)}{\mu_{m}(T)}=
 \frac{\mathcal H^m(S)}{\mathcal H^m(T)}
 =\frac{\lambda_{m}(S)}{\lambda_{m}(T)}\,.
 $$
 In conclusion, 
 $$
 \mu_{m}(S)=\lambda_m(S)\frac{\mu_{m}(T)}{\lambda_{m}(T)}=\lambda_m(S)\,.
 $$
The proof is complete.
 \end{proof}
 
 \section{Acknowledgement}
 The author is grateful to professors Andrew Glass
 and  Rolando Magnanini
 for valuable discussions on
  $GL(n,\mathbb Z) \ltimes \mathbb Z^{n}$ 
 and  $\mathcal H^n$, respectively.

\bibliographystyle{plain}

\end{document}